\documentclass[12pt]{article}

\usepackage{amsmath,amsthm}
\usepackage{amssymb}
\usepackage{hyperref}

\allowdisplaybreaks

\begin{document}

\title{Continued fraction expansions of the generating functions of Bernoulli and related numbers
}  
\author{
Takao Komatsu\\
}

\date{
\small MR Subject Classifications: Primary 11A55; Secondary 11J70, 30B70, 11B34, 11B68, 11B75, 05A19 
}

\maketitle

\def\fl#1{\left\lfloor#1\right\rfloor}
\def\cl#1{\left\lceil#1\right\rceil}
\def\ang#1{\left\langle#1\right\rangle}
\def\stf#1#2{\left[#1\atop#2\right]} 
\def\sts#1#2{\left\{#1\atop#2\right\}}

\newtheorem{theorem}{Theorem}
\newtheorem{Prop}{Proposition}
\newtheorem{Cor}{Corollary}
\newtheorem{Lem}{Lemma}

\begin{abstract} 
We give continued fraction expansions of the generating functions of Bernoulli numbers, Cauchy numbers, Euler numbers, harmonic numbers, and their generalized or related numbers.  In particular, we focus on explicit forms of the convergents of these continued fraction expansions. Linear fractional transformations of such continued fractions are also discussed.  We show more continued fraction expansion for different numbers and types, in particular, on Cauchy numbers. 
\end{abstract}


\section{Continued fractions} 

Given an analytic function $f(x)$, several types of general continued fractions have been known and studied. J. H. Lambert expanded $\log(1-x)$, $\arctan x$ and $\tan x$ in continued fractions in 1768.  Some special types are known as $J$-fractions, $C$-fractions, $T$-fractions, $M$-fractions (see, e.g., \cite{JT,Perron,Wall}). Lando \cite{Lando} discuss other similar continued fractions from the aspect of Combinatorics and generating functions.  Loya \cite[Chapter 8]{Loya} exhibits intriguing ideas with many useful examples of continued fraction expansions of analytic functions and famous numbers.  
However, it does not seem that explicit forms of the convergents of these generalized continued fraction expansions have been studied in detail.  

Assume that $f(x)$ is expanded by the Taylor expansion as $f(x)=\sum_{i=0}^\infty f_i x^i$ for complex numbers $f_i$ ($i=0,1,2,\dots$).   
In this paper, we focus on the following continued fraction expansion, which is related to $T$-fractions.  
\begin{equation}  
f(x)=1-\cfrac{h_1 x}{g_1+h_1 x-\cfrac{g_1 h_2 x}{g_2+h_2 x-\cfrac{g_2 h_3 x}{g_3+h_3 x-\cfrac{g_3 h_4 x}{g_4+h_4 x-{\atop\ddots}}}}}\,. 
\label{cf:genth}
\end{equation}
(Each concrete case needs some modifications.)   Since $a_n(x)=g_n+h_n x$ and $b_n(x)=-g_{n-1}h_n x$ ($n\ge 1$) with $g_0=1$ and $a_0(x)=1$, the convergents $P_n(x)/Q_n(x)$ of the continued fraction expansion of (\ref{cf:genth}) satisfy the recurrence relations 
\begin{align*} 
P_n(x)&=(g_n+h_n x)P_{n-1}(x)-g_{n-1}h_n x P_{n-2}(x)\quad(n\ge 1),\\ 
Q_n(x)&=(g_n+h_n x)Q_{n-1}(x)-g_{n-1}h_n x Q_{n-2}(x)\quad(n\ge 1)
\end{align*}
with $P_{-1}(x)=1$, $Q_{-1}(x)=0$, $P_0(x)=Q_0(x)=1$. 
By using the induction on $n$, we can know the explicit forms of $P$'s and $Q$'s as 
\begin{equation}  
P_n=g_1\cdots g_n\quad\hbox{and}\quad Q_n=g_1\cdots g_n\sum_{j=0}^n\frac{h_1\cdots h_j}{g_1\cdots g_j}x^j\,,
\label{pq:genth} 
\end{equation}  
respectively.  For convenience, $g_1\cdots g_j=h_1\cdots h_j=1$ when $j=0$. 
By the approximation property,  
\begin{equation}
Q_n f(x)\sim P_n(x)\quad(n\to\infty) 
\label{approx-pq}
\end{equation} 
or as polynomials 
$$ 
\frac{P_n(x)}{Q_n(x)}\equiv f(x)\pmod{x^{n+1}}
$$ 
Namely, 
$$
f(x)=\lim_{n\to\infty}\frac{P_n}{Q_n}=\left(\sum_{j=0}^\infty\frac{h_1\cdots h_j}{g_1\cdots g_j}x^j\right)^{-1}\,.
$$ 
In other words, the analytic function 
$$
f(x):=\sum_{n=i}^\infty f_i x^i=\left(\sum_{j=0}^\infty\frac{h_1\cdots h_j}{g_1\cdots g_j}x^j\right)^{-1}
$$ 
can be expanded as a continued fraction of type (\ref{cf:genth}).  In addition, Taylor expansion of $P_n/Q_n$ matches the terms of $f(x)$, up to $x^n$. 

\noindent 
{\bf Example.}  
Since 
$$
\frac{z}{\arctan z}=\left(\sum_{i=0}^\infty\frac{(-1)^i z^{2 i}}{2 i+1}\right)^{-1}\,, 
$$ 
the convergents $P_m(x)/Q_m(x)$ ($m\ge 0$) with 
$$
P_m(x)=(2 m+1)!!\quad\hbox{and}\quad Q_m(x)=(2 m+1)!!\sum_{i=0}^m\frac{(-1)^i z^{2 i}}{2 i+1}
$$ 
yield the continued fraction expansion 
$$
\frac{z}{\arctan z}=1+\cfrac{z^2}{3-z^2+\cfrac{(3 z)^2}{5-3 z^2+\cfrac{(5 z)^2}{7-5 z^2+\cfrac{(7 z)^2}{9-7 z^2+{\atop\ddots}}}}}\,. 
$$ 
In other words, we put $g_n=2 n+1$ ($n\ge 0$) and $h_n=-g_{n-1}=-(2 n-1)$ ($n\ge 1$) and $x=z^2$ in (\ref{cf:genth}). 
Thus, we get a famous generalized continued fraction of $\arctan z$, which is believed to be due to Leonhard Euler:     
$$
\arctan z=\cfrac{z}{1+\cfrac{z^2}{3-z^2+\cfrac{(3 z)^2}{5-3 z^2+\cfrac{(5 z)^2}{7-5 z^2+\cfrac{(7 z)^2}{9-7 z^2+{\atop\ddots}}}}}}\,. 
$$ 
In fact, the continued fraction expansion 
$$
\arctan z=\cfrac{z}{1+\cfrac{z^2}{3+\cfrac{(2 z)^2}{5+\cfrac{(3 z)^2}{7+\cfrac{(4 z)^2}{9+{\atop\ddots}}}}}}
$$ 
is more well-known (\cite[(6.1.14)]{JT},\cite[(90.3)]{Wall}), and is due to C. F. Gauss.  However, the convergents of this continued fraction is more complicated.  This is another reason why we do not consider different types of continued fraction expansions other than (\ref{cf:genth}) in this paper.

\section{Bernoulli numbers}  

Bernoulli numbers $B_n$ are defined by  
\begin{equation}
\sum_{n=1}^\infty B_n\frac{x^n}{n!}\,.  
\label{def:ber} 
\end{equation}  
Many kinds of continued fraction expansions of the generating functions of Bernoulli numbers have been known and studied (see, e.g., \cite[Appendix]{AIK},\cite{Frame}).  However, those of generalized Bernoulli numbers seem to be few, though there exist several generalizations of the original Bernoulli numbers.  In particular, direct generalizations from the continued fraction expansions seem to be hard.  

{\it Hypergeometric degenerate Bernoulli numbers} $\beta_{N,n}(\lambda)$ are defined by 
\begin{equation}   
\left({}_2F_1(1,N-1/\lambda;N+1;-\lambda x)\right)^{-1}=\sum_{n=0}^\infty\beta_{N,n}(\lambda)\frac{x^n}{n!} 
\label{def:deghbn} 
\end{equation}  
where 
$$
{}_2F_1(a,b;c;z)=\sum_{n=0}^\infty\frac{(a)^{(n)}b^{(n)}}{(c)^{(n)}}\frac{z^n}{n!}
$$ 
is the Gauss hypergeometric function with the rising factorial $(x)^{(n)}=x(x+1)\cdots(x+n-1)$ ($n\ge 1$) and $(x)^{(0)}=1$.  Denote the generalized falling factorial by $(x|r)_n=x(x-r)\cdots\bigl(x-(n-1)r\bigr)$ ($n\ge 1$) with $(x|r)_0=1$. $(x)_n=(x|1)_n$ is the original falling factorial. 
Since 
$$
{}_2F_1(1,N-1/\lambda;N+1;-\lambda x)=\sum_{n=0}^\infty\frac{(1-N\lambda|\lambda)_n}{(N+n)_n}x^n\,, 
$$ 
we consider the convergents $P_m(x)^\ast/Q_m(x)^\ast$ with 
$$
P_m(x)^\ast=1\quad\hbox{and}\quad Q_m(x)^\ast=\sum_{n=0}^m\frac{(1-N\lambda|\lambda)_n}{(N+n)_n}x^n\,.
$$ 
In order to treat with integer coefficients for $Q$'s, we can also consider the convergents $P_m(x)/Q_m(x)$, where 
$$
P_m(x)=(N+m)_m\quad\hbox{and}\quad Q_m(x)=(N+m)_m\sum_{n=0}^m\frac{(1-N\lambda|\lambda)_n}{(N+n)_n}x^n\,.
$$ 
Then 
\begin{align*}  
\frac{P_0(x)}{Q_0(x)}&=\frac{1}{1}=1,\quad 
\frac{P_1(x)}{Q_1(x)}=\frac{N+1}{N+1+(1-N\lambda)x}=1-\frac{(1-N\lambda)x}{N+1+(1-N\lambda)x},\\
\frac{P_2(x)}{Q_2(x)}&=\frac{(N+1)(N+2)}{(N+1)(N+2)+(N+2)(1-N\lambda)x+(1-N\lambda)_2 x^2}\\
&=1-\cfrac{(1-N\lambda)x}{(N+1)+(1-N\lambda)x-\cfrac{(N+1)\big(1-(N+1)\lambda\bigr)x}{(N+2)+\bigl(1-(N+1)\lambda\bigr)x}}
\end{align*} 
and $P_n(x)$ and $Q_n(x)$ ($n\ge 2$) satisfy the recurrence relations 
\begin{align*}  
&P_n(x)=\left(N+n+\bigl(1-(N+n-1)\lambda\bigr)x\right)P_{n-1}(x)\\
&\qquad\qquad\qquad\qquad\qquad -(N+n-1)\bigl(1-(N+n-1)\lambda\bigr)x P_{n-2}(x)\,,\\ 
&Q_n(x)=\left(N+n+\bigl(1-(N+n-1)\lambda\bigr)x\right)Q_{n-1}(x)\\
&\qquad\qquad\qquad\qquad\qquad -(N+n-1)\bigl(1-(N+n-1)\lambda\bigr)x Q_{n-2}(x)\,. 
\end{align*} 
By the recurrence relations, we know that $a_n(x)=N+n+\bigl(1-(N+n-1)\lambda\bigr)x$ ($n\ge 1$) and $b_n(x)=-(N+n-1)\bigl(1-(N+n-1)\lambda\bigr)x$ ($n\ge 2$) with $a_0(x)=1$ and $b_1(x)=-(1-N\lambda)x$. 
Therefore, we have the following continued fraction expansion of the generating function of the hypergeometric degenerate Bernoulli numbers.  

\begin{theorem}  
{\scriptsize 
\begin{align*}
&\sum_{n=0}^\infty\beta_{N,n}(\lambda)\frac{x^n}{n!}\\ 
&=1-\cfrac{(1-N\lambda)x}{N+1+(1-N\lambda)x-\cfrac{(N+1)\bigl(1-(N+1)\lambda\bigr)x}{N+2+\bigl(1-(N+1)\lambda\bigr)x-\cfrac{(N+2)\bigl(1-(N+2)\lambda\bigr)x}{N+3+\bigl(1-(N+2)\lambda\bigr)x-{\atop\ddots}}}}\,. 
\end{align*} 
}  
\label{cf:deghbn} 
\end{theorem}   

When $\lambda\to 0$ in Theorem \ref{cf:deghbn}, we get a continued fraction expansion of the generating function of hypergeometric Bernoulli numbers.  

\begin{Cor}  
$$
\sum_{n=0}^\infty B_{N,n}\frac{x^n}{n!}
=1-\cfrac{x}{N+1+x-\cfrac{(N+1)x}{N+2+x-\cfrac{(N+2)x}{N+3+x-{\atop\ddots}}}}\,. 
$$ 
\label{cf:hbn} 
\end{Cor}   

When $\lambda\to 0$ and $N=1$ in Theorem \ref{cf:deghbn}, we get a continued fraction expansion of the generating function of the original Bernoulli numbers defined in (\ref{def:ber}).  

\begin{Cor}  
$$
\sum_{n=0}^\infty B_{n}\frac{x^n}{n!}
=1-\cfrac{x}{2+x-\cfrac{2 x}{3+x-\cfrac{3 x}{4+x-\cfrac{4 x}{5+x-{\atop\ddots}}}}}\,. 
$$ 
\label{cf:bn} 
\end{Cor}   

We shall see the approximation property concerning Corollary \ref{cf:bn}.  
For example, for $n=4,5,6$ the convergents have Taylor expansions as 
\begin{align*}  
\frac{P_4(x)}{Q_4(x)}&=1-\frac{x}{2}+\frac{x^2}{12}-\frac{x^4}{720}+\frac{x^5}{720}-\frac{x^6}{864}+\frac{7 x^7}{17280}-\frac{29 x^8}{518400}-\cdots\,,\\
\frac{P_5(x)}{Q_5(x)}&=1-\frac{x}{2}+\frac{x^2}{12}-\frac{x^4}{720}+\frac{x^6}{4320}-\frac{x^7}{5760}+\frac{31 x^8}{518400}-\cdots\,,\\
\frac{P_6(x)}{Q_6(x)}&=1-\frac{x}{2}+\frac{x^2}{12}-\frac{x^4}{720}+\frac{x^6}{30240}+\frac{x^7}{40320}-\frac{83 x^8}{3628800}+\cdots\,. 
\end{align*} 
We see that 
$$
\sum_{n=0}^\infty B_n\frac{x^n}{n!}=1-\frac{x}{2}+\frac{x^2}{12}-\frac{x^4}{720}+\frac{x^6}{30240}-\frac{x^8}{1209600}+\cdots
$$ 

Thus, as expected, the coefficients of the Taylor expansion of $P_n(x)/Q_n(x)$ match those of the generating function of Bernoulli numbers up to the term of $x^n$.

Finally, when $\lambda\to 0$ and $N=0$ in Theorem \ref{cf:deghbn}, we get a continued fraction expansion of $e^{-x}$.  

\begin{Cor}  
$$
e^{-x}
=1-\cfrac{x}{1+x-\cfrac{x}{2+x-\cfrac{2 x}{3+x-\cfrac{3 x}{4+x-{\atop\ddots}}}}}\,. 
$$ 
\label{cf:e} 
\end{Cor}

\section{Cauchy numbers}  

Cauchy numbers are defined by 
\begin{equation}
\frac{x}{\log(1+x)}=\sum_{n=0}^\infty c_n\frac{x^n}{n!}\,. 
\label{def:cau} 
\end{equation}  
In \cite{DK}, new identities are found by using the continued fraction  
$$   
\frac{x}{\log(1+x)}=1+\cfrac{1^2 x}{2+\cfrac{1^2 x}{3+\cfrac{2^2 x}{4+\cfrac{2^2 x}{5+\cfrac{3^2 x}{6+{\atop\ddots}}}}}}
$$ 
(see, e.g., \cite[(90.1)]{Wall}).  
In this section, we find a different type of continued fraction expansion related to Cauchy numbers. 

Define the {\it hypergeometric degenerate Cauchy numbers} $\gamma_{N,n}(\lambda)$ by the generating function 
\begin{align}  
\frac{1}{{}_2 F_1(1,N-\lambda;N+1;-x)}
&=\frac{(\lambda-1)_{N-1}x^N/N!}{\frac{(1+x)^{\lambda}-1}{\lambda}-\sum_{n=1}^{N-1}(\lambda-1)_{n-1}x^n/n!}\notag\\
&:=\sum_{n=0}^\infty\gamma_{N,n}(\lambda)\frac{x^n}{n!}\,.  
\label{def:hgcn} 
\end{align}
When $\lambda\to 0$ in (\ref{def:hgcn}), $c_{N,n}=\gamma_{N,n}(0)$ are the hypergeometric Cauchy numbers (\cite{Ko3}).  When $N=1$ in (\ref{def:hgcn}), $\gamma_n(\lambda)=\gamma_{1,n}(\lambda)$ are the degenerate Cauchy numbers (see, e.g., \cite[(2.13)]{CH}). When $\lambda\to 0$ and $N=1$ in (\ref{def:hgcn}), $c_n=\gamma_{1,n}(0)$ are the classical Cauchy numbers in (\ref{def:cau}). $b_n=c_n/n!$ are often called Bernoulli numbers of the second kind.

Since 
$$ 
{}_2 F_1(1,N-\lambda;N+1;-x)=\sum_{n=0}^\infty\frac{(N-\lambda)^{(n)}}{(N+1)^{(n)}}(-t)^n
=\sum_{n=0}^\infty\frac{(\lambda-N)_n N!}{(N+n)!}x^n\,, 
$$ 
we consider the convergents $P_m(x)/Q_m(x)$ as 
$$
P_m(x)=\frac{(N+m)!}{N!}\quad\hbox{and}\quad Q_m(x)=(N+m)!\sum_{i=0}^m\frac{(\lambda-N)_i}{(N+i)!}x^i\,. 
$$ 
Then, we have 
\begin{align*}  
\frac{P_0(x)}{Q_0(x)}&=\frac{1}{1}=1,\quad 
\frac{P_1(x)}{Q_1(x)}=\frac{N+1}{N+1+(\lambda-N)x}=1-\frac{(\lambda-N)x}{N+1+(\lambda-N)x},\\
\frac{P_2(x)}{Q_2(x)}&=\frac{(N+1)(N+2)}{(N+1)(N+2)+(N+2)(\lambda-N)x+(\lambda-N)(\lambda-N-1)x^2}\\
&=1-\cfrac{(\lambda-N)x}{N+1+(\lambda-N)x-\cfrac{(N+1)(\lambda-N-1)x}{N+2+(\lambda-N-1)x}}
\end{align*}
and the convergents $P_n(x)/Q_n(x)$ ($n\ge 2$) satisfy the recurrence relations 
\begin{align*}  
&P_n(x)=\bigl(N+n+(\lambda-N-n+1)x\bigr)P_{n-1}(x)\\
&\qquad\qquad\qquad\qquad\qquad -(N+n-1)(\lambda-N-n+1)x P_{n-2}(x)\,,\\ 
&Q_n(x)=\bigl(N+n+(\lambda-N-n+1)x\bigr)Q_{n-1}(x)\\
&\qquad\qquad\qquad\qquad\qquad -(N+n-1)(\lambda-N-n+1)x Q_{n-2}(x)\,. 
\end{align*}
Since $a_n(x)=N+n+(\lambda-N-n+1)x$ ($n\ge 1$) and $b_n(x)=-(N+n-1)(\lambda-N-n+1)x$  ($n\ge 2$) with $a_0(x)=1$ and $b_1(x)=-(\lambda-N)x$, we have the following continued fraction expansion of the generating function of hypergeometric degenerate Cauchy numbers.  

\begin{theorem}  
{\small 
\begin{align*}
&\sum_{n=0}^\infty\gamma_{N,n}(\lambda)\frac{x^n}{n!}\\ 
&=1-\cfrac{(\lambda-N)x}{N+1+(\lambda-N)x-\cfrac{(N+1)(\lambda-N-1)x}{N+2+(\lambda-N-1)x-\cfrac{(N+2)(\lambda-N-2)x}{N+3+(\lambda-N-2)x-{\atop\ddots}}}}\,. 
\end{align*} 
}  
\label{cf:hdcn}
\end{theorem}

When $\lambda\to 0$ in Theorem \ref{cf:hdcn}, we have a continued fraction expansion of hypergeometric Cauchy numbers.

\begin{Cor}  
$$
\sum_{n=0}^\infty c_{N,n}\frac{x^n}{n!}
=1+\cfrac{N x}{N+1-N x+\cfrac{(N+1)^2 x}{N+2-(N+1)x+\cfrac{(N+2)^2 x}{N+3-(N+2)x-{\atop\ddots}}}}\,. 
$$ 
\label{cf:hcn} 
\end{Cor} 

When $\lambda\to 0$ and $N=1$ in Theorem \ref{cf:hdcn}, we have a continued fraction expansion of the original Cauchy numbers. 
\begin{align}
\sum_{n=0}^\infty c_{n}\frac{x^n}{n!}&=\frac{x}{\log(1+x)}\notag\\
&=1+\cfrac{x}{2-x+\cfrac{2^2 x}{3-2 x+\cfrac{3^2 x}{4-3 x+{\atop\ddots}}}}
\label{cf:cn}  
\end{align} 
(\cite[Chapter 8]{Loya}).

\section{Euler numbers}@@

{\it Hypergeometric Euler numbers} $E_{N,n}$ (\cite{KZ}) are defined by 
\begin{equation}
\frac{1}{{}_1 F_2(1;N+1,(2 N+1)/2;x^2/4)}=\sum_{n=0}^\infty E_{N,n}\frac{x^n}{n!}\,,
\label{def1:hypergeuler}
\end{equation}
where ${}_1 F_2(a;b,c;z)$ is the hypergeometric function defined by
$$
{}_1 F_2(a;b,c;z)=\sum_{n=0}^\infty\frac{(a)^{(n)}}{(b)^{(n)}(c)^{(n)}}\frac{z^n}{n!}\,.
$$
When $N=0$, then $E_n=E_{0,n}$ are classical Euler numbers, defined by 
$$
\frac{1}{\cosh x}=\sum_{n=0}^\infty E_n\frac{x^n}{n!}\,. 
$$  
Since by (\ref{def1:hypergeuler}) 
$$
{}_1 F_2(1;N+1,(2 N+1)/2;x^2/4)=\sum_{n=0}^\infty\frac{(2 N)!x^{2 n}}{(2 N+2 n)!}\,,
$$ 
consider the convergents $P_m(x)/Q_m(x)$ with  
$$
P_m(x)=\frac{(2 N+2 m)!}{(2 N)!}\quad\hbox{and}\quad Q_m(x)=(2 N+2 m)!\sum_{n=0}^m\frac{x^{2 n}}{(2 N+2 n)!}\,. 
$$ 
Now, 
\begin{align*}  
\frac{P_0(x)}{Q_0(x)}&=\frac{1}{1}=1,\\ 
\frac{P_1(x)}{Q_1(x)}&=\frac{(2 N+1)(2 N+2)}{(2 N+1)(2 N+2)+x^2}=1-\frac{x^2}{(2 N+1)(2 N+2)+x^2},\\ 
\frac{P_2(x)}{Q_2(x)}&=\frac{(2 N+1)(2 N+2)(2 N+3)(2 N+4)}{(2 N+1)(2 N+2)(2 N+3)(2 N+4)+(2 N+3)(2 N+4)x^2+x^4}\\ 
&=1-\cfrac{x^2}{(2 N+1)(2 N+2)+x^2-\cfrac{(2 N+1)(2 N+2)x^2}{(2 N+3)(2 N+4)+x^2}}
\end{align*}
and $P_n(x)/Q_n(x)$ ($n\ge 2$) satisfy the recurrence relations 
\begin{align*}  
P_n(x)&=\bigl((2 N+2 n-1)(2 N+2 n)+x^2\bigr)P_{n-1}(x)\\
&\qquad\qquad\qquad\qquad\qquad -(2 N+2 n-3)(2 N+2 n-2)x^2 P_{n-2}(x),\\
Q_n(x)&=\bigl((2 N+2 n-1)(2 N+2 n)+x^2\bigr)Q_{n-1}(x)\\
&\qquad\qquad\qquad\qquad\qquad -(2 N+2 n-3)(2 N+2 n-2)x^2 Q_{n-2}(x)\,. 
\end{align*}
Since $a_n(x)=(2 N+2 n-1)(2 N+2 n)+x^2$ ($n\ge 1$) and $b_n(x)=-(2 N+2 n-3)(2 N+2 n-2)x^2$ ($n\ge 2$) with $a_0(x)=1$ and $b_1(x)=-x^2$, we have the following continued fraction expansion of hypergeometric Euler numbers.

\begin{theorem}  
{\scriptsize 
\begin{align*} 
&\sum_{n=0}^\infty E_{N,n}\frac{x^n}{n!}\\
&=1-\cfrac{x^2}{(2 N+1)(2 N+2)+x^2-\cfrac{(2 N+1)(2 N+2)x^2}{(2 N+3)(2 N+4)+x^2-\cfrac{(2 N+3)(2 N+4)x^2}{(2 N+5)(2 N+6)+x^2-{\atop\ddots}}}}\,. 
\end{align*}
} 
\label{cf:hen} 
\end{theorem}

When $N=0$ in Theorem \ref{cf:hen}, we get a continued fraction expansion of the classical Euler numbers. 

\begin{Cor}  
\begin{align*} 
\sum_{n=0}^\infty E_{n}\frac{x^n}{n!}&=\frac{1}{\cosh x}\\
&=1-\cfrac{x^2}{1\cdot 2+x^2-\cfrac{1\cdot 2 x^2}{3\cdot 4+x^2-\cfrac{3\cdot 4 x^2}{5\cdot 6+x^2-{\atop\ddots}}}}\,. 
\end{align*}
\end{Cor}

{\it Hypergeometric Euler numbers of the second kind} $\widehat E_{N,n}$ (\cite{Ko10,KZ}) are defined by 
$$ 
\frac{1}{{}_1 F_2(1;N+1,(2 N+3)/2;x^2/4)}=\sum_{n=0}^\infty\widehat E_{N,n}\frac{x^n}{n!}\,.
$$
When $N=0$, $\widehat E_n=\widehat E_{0,n}$ are Euler numbers of the second kind or complementary Euler numbers, defined by 
$$
\frac{x}{\sinh x}=\sum_{n=0}^\infty\widehat E_n\frac{x^n}{n!} 
$$   
(\cite{Ko10}).  
Similarly, we have the following continued fraction expansions.  

\begin{theorem}  
{\scriptsize 
\begin{align*} 
&\sum_{n=0}^\infty\widehat E_{N,n}\frac{x^n}{n!}\\
&=1-\cfrac{x^2}{(2 N+2)(2 N+3)+x^2-\cfrac{(2 N+2)(2 N+3)x^2}{(2 N+4)(2 N+5)+x^2-\cfrac{(2 N+4)(2 N+5)x^2}{(2 N+6)(2 N+7)+x^2-{\atop\ddots}}}}\,. 
\end{align*}
} 
\label{cf:hen2} 
\end{theorem}  

\begin{Cor}  
\begin{align*} 
\sum_{n=0}^\infty\widehat E_{n}\frac{x^n}{n!}&=\frac{x}{\sinh x}\\
&=1-\cfrac{x^2}{2\cdot 3+x^2-\cfrac{2\cdot 3 x^2}{4\cdot 5+x^2-\cfrac{4\cdot 5 x^2}{6\cdot 7+x^2-{\atop\ddots}}}}\,. 
\end{align*}
\end{Cor}

\section{Harmonic numbers}  

There are several generalizations of harmonic numbers $H_n$, defined by 
$$
H_n=\sum_{k=1}^n\frac{1}{k}\quad(n\ge 1)\quad\hbox{with}\quad H_0=0\,. 
$$  
In \cite{Wang}, the generalized harmonic numbers of order $m$ ($m\ge 1$) are defined by 
\begin{equation}
h_n^{(m)}(a,b)=\sum_{k=1}^n\frac{1}{\bigl((k-1)a+b\bigr)^m}\quad(n\ge 1)\quad\hbox{with}\quad h_0^{(m)}(a,b)=0\,,  
\label{def:hno}
\end{equation}
where $a$ and $b$ are positive real numbers. 
When $a=b=1$, $H_n^{(m)}=h_n^{(m)}(1,1)$ are the $m$-order harmonic numbers.  
When $m=a=b=1$, $H_n=h_n^{(1)}(1,1)$ are the original harmonic numbers.  

Let $(x|r)^{(n)}=x(x+r)(x+2 r)\cdots\bigl(x+(n-1)r\bigr)$ ($n\ge 1$) be the generalized rising factorial with $(x|r)^{(0)}=1$.  When $r=1$, $(x)^{(n)}=(x|1)^{(n)}$ is the original rising factorial.  
From the definition in (\ref{def:hno}), we have 
$$
\sum_{n=1}^\infty h_n^{(m)}(a,b)x^n=\frac{1}{1-x}\sum_{k=1}^\infty\frac{x^k}{\bigl((k-1)a+b\bigr)^m}\,.
$$
Thus,  
$$
P_M'(x)=\sum_{k=1}^M\frac{x^k}{\bigl((k-1)a+b\bigr)^m},\quad Q_M'(x)=1-x
$$ 
or 
$$
P_M(x)=\bigl((b|a)^{(M)}\bigr)^m\sum_{k=1}^M\frac{x^k}{\bigl((k-1)a+b\bigr)^m},\quad Q_M(x)=\bigl((b|a)^{(M)}\bigr)^m(1-x)
$$ 
yield that 
$$
Q_M'(x)\left(\sum_{n=1}^\infty h_n^{(m)}(a,b)x^n\right)\sim P_M'(x)\quad(M\to\infty)
$$ 
or 
$$
Q_M(x)\left(\sum_{n=1}^\infty h_n^{(m)}(a,b)x^n\right)\sim P_M(x)\quad(M\to\infty)\,. 
$$ 

Now, 
{\small 
\begin{align*}
\frac{P_0(x)}{Q_0(x)}&=\frac{0}{1-x}=0,\quad \frac{P_1(x)}{Q_1(x)}=\frac{x}{b^m(1-x)}\,,\\ 
\frac{P_2(x)}{Q_2(x)}&=\frac{(a+b)^m x+b^m x^2}{b^m(a+b)^m(1-x)}\\
&=\cfrac{x}{b^m(1-x)-\cfrac{b^{2 m}x(1-x)}{(a+b)^m+b^m x}}\,,\\
\frac{P_3(x)}{Q_3(x)}&=\frac{(a+b)^m(2 a+b)^m x+b^m(2 a+b)^m x^2+b^m(a+b)^m x^3}{b^m(a+b)^m(2 a+b)^m(1-x)}\\
&=\cfrac{x}{b^m(1-x)-\cfrac{b^{2 m}x(1-x)}{(a+b)^m+b^m x-\cfrac{(a+b)^{2 m}x}{(2 a+b)^m+(a+b)^m x}}}\,,
\end{align*} 
} 
and $P_n(x)$ and $Q_n(x)$ ($n\ge 3$) satisfy the recurrence relations 
{\small 
\begin{align*}  
&P_n(x)=\left(\bigl((n-1)a+b\bigr)^m+\bigl((n-2)a+b\bigr)^m x\right)P_{n-1}(x)\\
&\qquad\qquad\qquad\qquad\qquad\qquad\qquad\qquad -\bigl((n-2)a+b\bigr)^{2 m}x P_{n-2}(x)\,,\\ 
&Q_n(x)=\left(\bigl((n-1)a+b\bigr)^m+\bigl((n-2)a+b\bigr)^m x\right)Q_{n-1}(x)\\
&\qquad\qquad\qquad\qquad\qquad\qquad\qquad\qquad -\bigl((n-2)a+b\bigr)^{2 m}x Q_{n-2}(x)\,.   
\end{align*} 
} 
It is proved by induction. The right-hand side of $P$'s is equal to 
\begin{align*}  
&\bigl((n-1)a+b\bigr)^m\bigl((b|a)^{(n-1)}\bigr)^m\sum_{k=1}^{n-1}\frac{x^k}{\bigl((k-1)a+b\bigr)^m}\\
&\quad+\bigl((n-2)a+b\bigr)^m x\bigl((b|a)^{(n-1)}\bigr)^m\sum_{k=1}^{n-1}\frac{x^k}{\bigl((k-1)a+b\bigr)^m}\\
&\quad -\bigl((n-2)a+b\bigr)^{2 m}x\bigl((b|a)^{(n-2)}\bigr)^m\sum_{k=1}^{n-2}\frac{x^k}{\bigl((k-1)a+b\bigr)^m}\\
&=\bigl((b|a)^{(n-1)}\bigr)^m\sum_{k=1}^{n-1}\frac{x^k}{\bigl((k-1)a+b\bigr)^m}\\
&\quad +\bigl((n-2)a+b\bigr)^m x\bigl((b|a)^{(n-1)}\bigr)^m\frac{x^{n-1}}{\bigl((n-2)a+b\bigr)^m}\\
&=\bigl((b|a)^{(n-1)}\bigr)^m\sum_{k=1}^{n}\frac{x^k}{\bigl((k-1)a+b\bigr)^m}=Q_n(x)\,. 
\end{align*} 
The relation for $Q$'s is trivially checked.  
 
Since 
$a_n(x)=\bigl((n-1)a+b\bigr)^m+\bigl((n-2)a+b\bigr)^m x$ ($n\ge 2$) and $b_n(x)=-\bigl((n-2)a+b\bigr)^{2 m}x$ ($n\ge 3$) 
with $a_0(x)=0$, $a_1(x)=b^m(1-x)$, $b_1(x)=x$, $b_2(x)=-b^{2 m}x(1-x)$,  
we have the following continued fraction expansion.  

\begin{theorem} 
{\scriptsize 
\begin{multline*} 
\sum_{n=1}^\infty h_n^{(m)}(a,b)x^n\\ 
=\cfrac{x}{b^m(1-x)-\cfrac{b^{2 m}x(1-x)}{(a+b)^m+b^m x-\cfrac{(a+b)^{2 m}x}{(2 a+b)^m+(a+b)^m x-\cfrac{(2 a+b)^{2 m}x}{(3 a+b)^m+(2 a+b)^m x-{\atop\ddots}}}}}\,.
\end{multline*}
}
\label{cf:ghn1} 
\end{theorem}

When $a=b=1$ in Theorem \ref{cf:ghn1}, we have a continued fraction expansion concerning the $m$-order harmonic numbers.  

\begin{Cor} 
$$ 
\sum_{n=0}^\infty H_n^{(m)}x^n 
=\cfrac{x}{(1-x)-\cfrac{x(1-x)}{2^m+x-\cfrac{2^{2 m}x}{3^m+2^m x-\cfrac{3^{2 m}x}{4^m+3^m x-{\atop\ddots}}}}}\,.
$$ 
\label{cf:ghn2}  
\end{Cor}  

When $m=a=b=1$ in Theorem \ref{cf:ghn1}, we have a continued fraction expansion concerning the original harmonic numbers.  

\begin{Cor} 
$$  
\sum_{n=0}^\infty H_n x^n 
=\cfrac{x}{(1-x)-\cfrac{x(1-x)}{2+x-\cfrac{2^{2}x}{3+2 x-\cfrac{3^{2}x}{4+3 x-{\atop\ddots}}}}}\,. 
$$ 
\label{cf:hn3}  
\end{Cor}

\section{Functions associated with the Riemann zeta function}  

Let $\mu(n)$ be M\"obius function.  From the property 
$$
\sum_{d|n}\mu(d)=\begin{cases}  
1&\text{if $n=1$},\\ 
0&\text{if $n\ge 2$} 
\end{cases}\,,
$$ 
the Dirichlet series that generates the M\"obius function is the multiplicative inverse of the Riemann zeta function:  
\begin{equation}  
\sum_{n=1}^\infty\frac{\mu(n)}{n^s}=\frac{1}{\zeta(s)}\,,
\label{gf:mobius}
\end{equation}  
where $s$ is a complex number with real part larger than $1$.  
Then we consider the convergents $P_m(x)/Q_m(x)$ as 
$$
P_m(x)=(m!)^s,\quad Q_m(x)=(m!)^s\sum_{k=1}^m\frac{x^k}{k^s}\quad(m\ge 1)\,. 
$$ 
Now, 
\begin{align*} 
&\frac{P_1(x)}{Q_1(x)}=\frac{1}{x},\quad 
\frac{P_2(x)}{Q_2(x)}=\frac{2^s}{2^s x+x^2}=\cfrac{1}{x+\frac{x^2}{2^s}},\\
&\frac{P_3(x)}{Q_3(x)}=\frac{6^n}{6^s x+3^s x^2+2^s x^3}=\cfrac{1}{x+\cfrac{x^2}{2^s-\cfrac{4^s x}{3^s+2^s x}}}
\end{align*}
and 
$P_n(x)$ and $Q_n(x)$ ($n\ge 3$) satisfy the recurrence relations 
{\small 
\begin{align*}  
&P_n(x)=\left(n^s+(n-1)^s x\right)P_{n-1}(x)-(n-1)^{2 s}x P_{n-2}(x)\,,\\ 
&Q_n(x)=\left(n^s+(n-1)^s x\right)Q_{n-1}(x)-(n-1)^{2 s}x Q_{n-2}(x)\,.   
\end{align*} 
Since $a_n(x)=$ and $b_n(x)=$ with $a_1(x)=x$, $a_2(x)=2^s$, $b_1(x)=1$ and $b_2(x)=x^2$, we have the continued fraction expansion 
$$ 
\lim_{n\to\infty}\frac{P_n(x)}{Q_n(x)}
=\cfrac{1}{x+\cfrac{x^2}{2^s-\cfrac{2^{2 s}x}{3^s+2^s x-\cfrac{3^{2 s} x}{4^s+3^s x-\cfrac{4^{2 s} x}{5^s+4^s x-{\atop\ddots}}}}}}\,.
$$  
Notice that 
$$
\lim_{n\to\infty}\frac{P_n(x)}{Q_n(x)}\ne\sum_{n=1}^\infty\frac{\mu(n)}{n^s}x^n\,. 
$$ 

Putting $x=1$, we have a continued fraction expansion of the Dirichlet series of M\"obius function.  

\begin{theorem}  
$$
\sum_{n=1}^\infty\frac{\mu(n)}{n^s}=\cfrac{1}{1+\cfrac{1}{2^s-\cfrac{2^{2 s}}{3^s+2^s-\cfrac{3^{2 s}}{4^s+3^s-\cfrac{4^{2 s}}{5^s+4^s-{\atop\ddots}}}}}}\,.
$$ 
\end{theorem}

For example, for $s=7$
$$
\frac{P_5(1)}{Q_5(1)}=0.9917254568069276497590711416
$$ 
but 
$$
\sum_{n=1}^\infty\frac{\mu(n)}{n^5}=\frac{1}{\zeta(7)}=0.9917198558384443104281859315\,. 
$$

\section{Transforms of continued fractions}  

Raney \cite{Raney} established the method to yield the simple continued fraction expansions of $\displaystyle \beta(x)=\frac{a x+b}{c x+d}$ with $a d-b c\ne 0$ from the simple continued fraction expansions 
$$
\alpha=[a_0;a_1,a_2,a_3,\dots]:=a_0+\cfrac{1}{a_1+\cfrac{1}{a_2+\cfrac{1}{a_3+{\atop\ddots}}}}\,, 
$$   
where $a_0$ is an integer and $a_1,a_2,\dots$ are positive integers. 
In this section,  we shall consider the linear fractional transformation $\beta(x)$ 
for the continued fraction expansion in (\ref{cf:genth}).  
Since the convergents $P_n(x)/Q_n(x)$ ($n\ge 0$) are given in (\ref{pq:genth}) 
in our case,  
we have 
\begin{align*}
\frac{\tilde P_n(x)}{\tilde Q_n(x)}&:=\beta\left(\frac{P_n(x)}{Q_n(x)}\right)\\
&=\dfrac{g_1\cdots g_n\left(\bigl(a+b\bigr)+b\sum_{j=1}^n\frac{h_1\cdots h_j}{g_1\cdots g_j}x^j\right)}{g_1\cdots g_n\left(\bigl(c+d\bigr)+d\sum_{j=1}^n\frac{h_1\cdots h_j}{g_1\cdots g_j}x^j\right)}\,. 
\end{align*}
Since 
\begin{align*}  
\frac{\tilde P_0(x)}{\tilde Q_0(x)}&=\frac{a+b}{c+d}\,,\\ 
\frac{\tilde P_1(x)}{\tilde Q_1(x)}&=\frac{(a+b)g_1+b h_1 x}{(c+d)g_1+d h_1 x}
=\frac{a+b}{c+d}-\dfrac{\dfrac{a d-b c}{c+d}h_1 x}{(c+d)g_1+d h_1 x}\,,\\
\frac{\tilde P_2(x)}{\tilde Q_2(x)}&=\frac{(a+b)g_1 g_2+b h_1 g_2 x+b h_1 h_2 x^2}{(c+d)g_1 g_2+d h_1 g_2 x+d h_1 h_2 x^2}\\
&=\frac{a+b}{c+d}-\cfrac{\dfrac{a d-b c}{c+d}h_1 x}{(c+d)g_1+d h_1 x-\cfrac{(c+d)g_1 h_2 x}{g_2+h_2 x}}
\end{align*}
and the convergents $\tilde P_n(x)/{\tilde Q}_n(x)$ ($n\ge 3$) satisfy the recurrence relations 
\begin{align*} 
\tilde P_n(x)&=(g_n+h_n x)\tilde P_{n-1}(x)-g_{n-1}h_n x\tilde P_{n-2}(x),\\ 
\tilde Q_n(x)&=(g_n+h_n x)\tilde Q_{n-1}(x)-g_{n-1}h_n x\tilde Q_{n-2}(x)\,.
\end{align*} 
Since $a_n(x)=g_n+h_n x$ ($n\ge 2$) and $b_n(x)=-g_{n-1}h_n x$ ($n\ge 3$) with $a_0(x)=(a+b)/(c+d)$, $a_1(x)=(c+d)g_1+d h_1 x$, $b_1(x)=-(a d-b c)/(c+d) h_1 x$ and $b_2(x)=-(c+d)g_1 h_2 x$, we have the following continued fraction expansion.  

\begin{theorem}  
If $f(x)$ has a continued fraction expansion (\ref{cf:genth}), we have 
$$ 
\beta\big(f(x)\bigr)=\frac{a+b}{c+d}-\cfrac{\dfrac{a d-b c}{c+d}h_1 x}{(c+d)g_1+d h_1 x-\cfrac{(c+d)g_1 h_2 x}{g_2+h_2 x-\cfrac{g_2 h_3 x}{g_3+h_3 x-\cfrac{g_3 h_4 x}{g_4+h_4 x-{\atop\ddots}}}}}\,. 
$$ 
\label{cf:transform}
\end{theorem} 

\noindent
{\bf Examples.}  
Applying Theorem \ref{cf:transform} to the continued fraction expansions in Corollary \ref{cf:bn}, Corollary \ref{cf:e} and (\ref{cf:cn}), we obtain the following expansions.  
When $a=-1$ and $b=c=d=1$, we get 
\begin{align*}
\beta\left(\sum_{n=0}^\infty B_{n}\frac{x^n}{n!}\right)&=\frac{e^x-x-1}{e^x+x-1}\\
&=\cfrac{x}{4+x-\cfrac{4 x}{3+x-\cfrac{3 x}{4+x-\cfrac{4 x}{5+x-{\atop\ddots}}}}} 
\end{align*}
and 
\begin{align*}
\beta(e^{-x})&=\tanh(x/2)\\ 
&=\cfrac{x}{2+x-\cfrac{2 x}{2+x-\cfrac{2 x}{3+x-\cfrac{3 x}{4+x-\cfrac{4 x}{5+x-\cfrac{5 x}{6+x-{\atop\ddots}}}}}}}\,. 
\end{align*}
When $b=-1$ and $a=c=d=1$, we get 
\begin{align*}
\beta\left(\sum_{n=0}^\infty c_{n}\frac{x^n}{n!}\right)&=\frac{x-\log(1+x)}{x+\log(1+x)}\\
&=\cfrac{x}{4-x+\cfrac{2^3 x}{3-2 x+\cfrac{3^2 x}{4-3 x+\cfrac{4^2 x}{5-4 x+{\atop\ddots}}}}}\,.   
\end{align*}

\subsection{Linear fractional transformation}  

Theorem \ref{cf:transform} can be proved by tools by Raney to manipulate expansions.  This method is due to Hurwitz (see, e.g., \cite{Hurwitz}), Frame \cite{Frame1949}, and Kolden \cite{Kolden}, independently, 
 but popularized by van der Poorten (see, e.g., \cite{vdP}).  
For simplicity, put $P_n=P_n(x)$, $Q_n=Q_n(x)$, $a_n=a_n(x)$ and $b_n=b_n(x)$. 
Since the recurrence relation can be explained by matrices as 
$$ 
\begin{pmatrix}a_0&1\\1&0\end{pmatrix}\begin{pmatrix}a_1&1\\b_1&0\end{pmatrix}\cdots\begin{pmatrix}a_n&1\\b_n&0\end{pmatrix}=\begin{pmatrix}P_n&Q_n\\P_{n-1}&Q_{n-1}\end{pmatrix}\,,
$$ 
the linear fractional transformation $\beta(x)=(a x+b)/(c x+d)$ implies the multiplication of the matrix 
$$
\begin{pmatrix}a&b\\c&d\end{pmatrix}
$$ 
from the left.  Namely,  if 
$$
\beta\left(\frac{P_n}{Q_n}\right)=\frac{\tilde P_n}{\tilde Q_n}\,,
$$ 
then 
\begin{align*} 
\begin{pmatrix}\tilde P_n&\tilde Q_n\\\tilde P_{n-1}&\tilde Q_{n-1}\end{pmatrix}
&:=\begin{pmatrix}\tilde a_0&1\\1&0\end{pmatrix}\begin{pmatrix}\tilde a_1&1\\\tilde b_1&0\end{pmatrix}\cdots\begin{pmatrix}\tilde a_n&1\\\tilde b_n&0\end{pmatrix}\\ 
&=\begin{pmatrix}a&b\\c&d\end{pmatrix}\begin{pmatrix}a_0&1\\1&0\end{pmatrix}\begin{pmatrix}a_1&1\\b_1&0\end{pmatrix}\cdots\begin{pmatrix}a_n&1\\b_n&0\end{pmatrix}\,. 
\end{align*}  
Now, we see that 
\begin{align*}
&\begin{pmatrix}a&b\\c&d\end{pmatrix}\begin{pmatrix}1&1\\1&0\end{pmatrix}\begin{pmatrix}g_1+h_1 x&1\\-h_1 x&0\end{pmatrix}\begin{pmatrix}g_2+h_2 x&1\\-g_1 h_2 x&0\end{pmatrix}\\
&=\begin{pmatrix}(a+b)g_1 g_2+b h_1 g_2 x+b h_1 h_2 x^2&(a+b)g_1+b h_1 x\\(c+d)g_1 g_2+d h_1 g_2 x+d h_1 h_2 x^2&(c+d)g_1+d h_1 x\end{pmatrix}\\
&=\begin{pmatrix}\frac{a+b}{c+d}&1\\1&0\end{pmatrix}\begin{pmatrix}(c+d)g_1+d h_1 x&1\\-\frac{a d-b c}{c+d}h_1 x&0\end{pmatrix}\begin{pmatrix}g_2+h_2 x&1\\-(c+d) g_1 h_2 x&0\end{pmatrix}\,. 
\end{align*}  
Therefore, we obtain the continued fraction expansion in Theorem \ref{cf:transform}.  

There are different types of transformations.  
In \cite{Barry}, the $r$-th binomial transform $b_n=\sum_{k=0}^n r^{n-k}\binom{n}{k}a_k$ and some more variations are performed in the sense of Riordan arrays $(1/(1-r x),x/(1-r x))$ and so on.

\section{Continued fractions of the ordinary generating functions}  

A continued fraction expansion of the ordinary generating function of Bernoulli numbers is given by 
$$
\sum_{n=0}^\infty B_n x^n
=\cfrac{1}{1+\cfrac{x}{\dfrac{2}{1}-\cfrac{x}{3+\cfrac{2 x}{\dfrac{2}{2}-\cfrac{2 x}{5+\cfrac{3 x}{\dfrac{2}{3}-\cfrac{3 x}{7+\cfrac{4 x}{\dfrac{2}{4}-\cfrac{4 x}{9+\dfrac{5 x}{\ddots}}}}}}}}}}
$$ 
(\cite[A.5]{AIK}). 
However, any beautiful continued fraction expansion for Cauchy numbers has not been known yet.  Nevertheless, in order to satisfy the approximation property (\ref{approx-pq}), we have the following expansion.  

\begin{theorem}  
\begin{align*}  
\sum_{n=0}^\infty c_n x^n&=1+\frac{x}{2}-\frac{x^2}{6}+\frac{x^3}{4}-\frac{19}{30}x^4+\frac{9}{4}x^5-\frac{863}{84}x^6+\frac{1375}{24}x^7-\frac{33953}{90}x^8
+\cdots\\
&=1+\cfrac{x}{2+\cfrac{2 x}{3+\cfrac{7 x}{2+\cfrac{93 x}{35+\cfrac{2391 x}{31+\cfrac{172542 x}{2391+\cfrac{443242433 x}{57514+\cfrac{19157845135465 x}{100087001+{\atop\ddots}}}}}}}}}\,.
\end{align*}
\label{cf:ogf-cn} 
\end{theorem}

\section{Concluding remarks}   

In this paper, we deal with continued fraction in the aspects of convergents. Such techniques and ideas can be applied to more different types.  For example, we can have a more complicated continued fraction expansion than that in Corollary \ref{cf:hcn}.   
  
\begin{theorem}  
For $N\ge 1$,  
\begin{align*}  
&\sum_{n=0}^\infty c_{N,n}\frac{x^n}{n!}\\
&=1+\cfrac{N x}{N+1+\cfrac{x}{N+2+\cfrac{(N+1)^2 x}{N+3+\cfrac{2^2 x}{N+4+\cfrac{(N+2)^2 x}{N+5+\cfrac{3^2 x}{N+6+\cfrac{(N+3)^2 x}{N+7+{\atop\ddots}}}}}}}}\,. 
\end{align*} 
\label{th:hcn} 
\end{theorem}   

When $N=1$, this is a direct generalization of the continued fraction expansion 
\begin{align*}   
\frac{x}{\log(1+x)}&\sum_{n=0}^\infty c_{n}\frac{x^n}{n!}\\
&=1+\cfrac{1^2 x}{2+\cfrac{1^2 x}{3+\cfrac{2^2 x}{4+\cfrac{2^2 x}{5+\cfrac{3^2 x}{6+{\atop\ddots}}}}}}\,.    
\end{align*}  
(see, e.g., \cite[(90.1)]{Wall}).  

Similarly, the $n$-th convergent and the generating function of $c_n$ coincide up to the $n$-th term in their Taylor expansions.  However, the structure of the continued fraction expansion in Theorem \ref{th:hcn} is more complicated than that in Corollary \ref{cf:hcn}.  In addition, the corresponding Euler continued fractions in Theorem \ref{cf:hen} and Theorem \ref{cf:hen2} have not been found yet.  The details will be discussed in the following papers.

\end{document}